\documentclass[9pt]{article}
\usepackage{amsfonts,amssymb,amsmath}
\usepackage[T2A]{fontenc}  
\usepackage[cp1251]{inputenc}
\usepackage[english]{babel}
\usepackage{geometry}
\usepackage{color}     
\newtheorem{teor}{Theorem}[section]

\newtheorem{lemma}[teor]{Lemma}

\numberwithin{equation}{section}

\usepackage{tocloft}

\def\ds{\displaystyle}
\def\f{\frac}
\def\div{{\rm div}}
\def\C{\mathcal{C}}
\def\R{\mathbb{R}}

\def\K{\mathcal{K}}
\def\L{\mathcal{L}}

\def\per{{\rm per}}

\def\nab{\nabla}
\def\Div{{\rm div}}
\def\pa{\partial}
\def\l{\left}
\def\r{\right}

\def\ve{\varepsilon}
\def\e{\varepsilon}
\def\o{\omega}
\def\C0{C^\infty_0(\R^d)}

\def\ld{L^2(\R^d)}
\def\rd{\R^d}
\def\RR^d{\R^d}

\def\beq{\begin{equation}}
\def\eeq{\end{equation}}

\begin{document}
\title{On resolvent approximations   of elliptic differential
 operators with locally periodic coefficients  
}

\author{ S.\,E. Pastukhova}
\date{}
\maketitle


\begin{footnotesize} 
We study the asymptotic behaviour of the resolvents $(A_\e+1)^{-1}$
of 
elliptic
second-order differential 
operators $A_\e=-\Div\, a^\e(x)\,\nab$ in  $\R^d$ with rapidly oscillating 
coefficients, as the small parameter $\e$ tends to zero. The matrix $a^\e(x)=a(x,x/\e)$ has the two-scale structure: it depends on the fast variable  $x/\e$ and on the slow variable $x$, with periodicity only in the fast variable.
We provide a construction for the leading  terms in the “operator
asymptotics” of $(A_\e+1)^{-1}$
 in the sense of $L^2$-operator-norm convergence 
 with 
order $\e^2$ remainder estimates. We apply the 
modified  method of the
first approximation with the usage of the shift proposed  by V.V. Zhikov in
\textit{Dokl. Math.}, {\bf 72}:1 (2005).

\end{footnotesize}
\bigskip

\section{Introduction}
\textbf{1.1. $L^2$-estimate of homogenization error in periodic setting.}
This paper relates to  homogenization theory, more precisely, to its branch connected with operator-type estimates for the error of homogenization.
Homogenization studies heterogeneous media via corresponding differential equations and integral functionals (for introduction see, for example, the books \cite{BLP}, \cite{SP}, \cite{BP},  \cite{ZKO}).

An elliptic equation with periodic rapidly oscillating coefficients is one of  model examples in homogenization
theory. It may describe various physical processes, say, in 
small-period composites. 
Suppose that the composite has only two phases (one can imagine grains that are periodically distributed in a host  medium)  and the contrast between them is moderate. In the simplest case, the latter means that the physical characteristics of the  both phases are positive  constants which are distinct. 
Theory of homogenisation 
aims at characterising limiting, or "effective", properties of small-period composites. A typical problem here
is to study the asymptotic behaviour as $\e\to 0$ of the solutions $u^\e$ to the equations of the type
\beq\label{0.1}
\ds{u^\ve\in H^1(\RR^d),\quad A_\ve u^\ve+ u^\ve=f,\quad
f\in L^2(\RR^d),} \atop\ds{A_\ve =-\div\, a(x/\e)\nab},\eeq
where the matrix
$a(x)$ is symmetric, measurable, $Y$-periodic, $Y=[-1/2,1/2)^d$, and satisfies the condition of uniform ellipticity: $\lambda I\le a(\cdot)\le \lambda^{-1} I$ for some positive constant $\lambda\le 1$.
A well known result says that the limiting, or "effective", equation  is of the same type but much simpler, namely,
\beq\label{0.2}
\ds{u\in H^1(\RR^d),\quad A_0 u+ u=f, \quad
f\in L^2(\RR^d),} \atop\ds{A_0 =-\div\, a^0\nab},\eeq
 with the constant matrix $a^0>0$ which is found by solving an  auxiliary problem on the cell of periodicity $Y$. The closeness of the solutions to the  problems (\ref{0.1}) and (\ref{0.2}) may be expressed in different forms and can be proved by different approaches; besides, the estimate of this closeness is of great interest, in particular, for applications. 
  The strongest result in this direction states that the resolvents $(A_\ve+1)^{-1}$ and $(A_0+1)^{-1}$ 
are close to each other in the $L^2$-operator norm; moreover, the following estimate sharp in order holds:
\beq\label{0.3}
\|(A_\ve+1)^{-1}-(A_0+1)^{-1}\|_{L^2\to L^2}\le C\e
\eeq
with constant $C$ depending only the dimension $d$ and the ellipticity constant $\lambda$. In other words, 
the estimate 
(\ref{0.3}) means that the difference of the solutions to the problems (\ref{0.1}) and (\ref{0.2}) satisfies the estimate
\beq\label{0.4}
\|u^\e-u\|_{L^2}\le C\e\|f\|_{L^2}, \quad C=const(d,\lambda). 
\eeq

In the framework of more general results, the operator-type estimate (\ref{0.3}) was proved for the first time in \cite{BS} and a little bit later  also  in \cite{Zh1}, in the latter by another method, which is much simpler conceptually than that of  \cite{BS}. 
Before this result was established, $L^2$-estimates for the difference
$u^\e-u$ had been proved under more restrictive assumptions on regularity   of the coefficients  in the equation (\ref{0.1}) 
(say, they should be from the space $C^k$, $k\in \mathbb{N}$ is sufficiently large)
and
 with majorants 
that
contain the high order  Sobolev norms of the right-hand side function $f$. Namely, these were the error estimates of the form
\[
\|u^\e-u\|_{L^2}\le C\e, 
\]
where the constant $C$ depended on the   Sobolev norm  $\|f\|_{H^m}$
of the right-hand side function
and the norm $\|a\|_{C^k}$ of the coefficients in the equation. Passing from such kind estimates to the operator estimate (\ref{0.3}) for the difference of the resolvents is impossible.

Another type of the operator convergence (different from the uniform resolvent convergence  discussed above) that connects the family $A_\e$ with  the effective operator $A_0$ is the strong resolvent convergence which is also often considered in homogenization (see discussion in \S2).

\textbf{1.2. 
Further extensions.}
The result  (\ref{0.3}) may be extended in different directions. First of all, the question arises whether it is possible, under the same minimal assumptions on regularity, 
to obtain approximations for $(A_\ve+1)^{-1}$ in the same $L^2$-operator norm with remainder terms of $\e^2$ order. The positive answer was  given in \cite{Zh05} and \cite{BS05} where, in the framework of more general results, the following estimate was proved:
\beq\label{0.5}
\ds{
\|(A_\ve+1)^{-1}-(A_0+1)^{-1} -\e\,\mathcal{C}_\e\|_{L^2\to L^2}\le C\e^2,} \atop\ds{\mathcal{C}_\e=\K_\e+(\K_\e)^*,
}
\eeq
with the constant $C$ of the same type as in (\ref{0.3}). Here 
 $\K_\e$ is the correcting operator such that 
\beq\label{0.5c}
\K_\e\, f (x)= N(x/\e) \cdot\nab u(x), \quad \mbox{ where}\quad u=(A_0+1)^{-1}f,
\eeq
 $N(\cdot)$ is a solution of some auxiliary periodic problem.

The corrector term (\ref{0.5c})  is  well known and widely used in  homogenization theory. For example, this corrector is taken in the classical homogenization, under high regularity conditions, to construct the $H^1$-approximation of the solution to the equation  (\ref{0.1}). That is the function $u^1_\e(x)=u(x)+\e N(x/\e) \cdot\nab u(x)$
which enables the approximation of $u^\e$ with the estimate 
\beq\label{0.5h}
\|u^\e-u^1_\e\|_{H^1}\le C\e.
\eeq
For a long time this estimate was obtained with the constant $C$
depending on additional regularity characteristics of the function $f$ and the matrix  $a$.

Assuming our minimal regularity conditions, even the existence of $u^1_\e$ as an element of the space
$H^1(\RR^d)$ is under the question, and it seems from the first sight that   some additional regularity on data is necessary if we want to write any
 $H^1$-approximations and prove estimates for them.
But it turns out that,
in the scalar case that is now at hand, 
 no additional regularity on data is needed to obtain a more general result than (\ref{0.5h}), namely, the following operator 
estimate:
\beq\label{0.6}
\|(A_\ve+1)^{-1}-(A_0+1)^{-1} -\e\,\K_\e\|_{L^2\to H^1}\le c\e,\quad c=const(d,\lambda),
\eeq
where the operator $\K_\e$ is the same as in (\ref{0.5}) and (\ref{0.5c}). So the estimate (\ref{0.5h}) acquires the form
\[
\|u^\e-u^1_\e\|_{H^1}\le c\e\|f\|_{L^2},\quad c=const(d,\lambda),
\]
wherefrom the $L^2$-estimate (\ref{0.4}) follows as a simple corollary.

The operator estimate (\ref{0.6}) as well as its more general counterpart, with Steklov's smoothing in corrector that is needed for vector problems, were proved for the first time in \cite{ZhP05}. 

Another direction to extend the result (\ref{0.3}) is to go off pure periodicity in coefficients  of the operator $A_\e$, letting them to be, for example, locally periodic (see the precise definition in \S2) which depend 
on the fast variable  $x/\e$, as $\e\to 0$, and on the slow variable $x$, with periodicity only in the fast variable.
For this setup of the problem, the  error estimate  of the type (\ref{0.3})  and (\ref{0.6}) have been already proved in \cite{PT}. Here, \textit{we focus  on the estimate of the type (\ref{0.5})  for a locally periodic operator $A_\e$ which is not necessarily selfadjoint. } In this case, certainly, the structure of the corrector in the approximation of the resolvent $(A_\ve+1)^{-1}$ with order $\e^2$ remainder becomes more complicated than in (\ref{0.5}): there emerge several additional terms in the corrector relating to the locally periodicity and the
nonselfadjointness either. 

The main result is formulated in Theorem \ref{Th5.1} and is proved in \S6. Preliminary lemmas 
that play a key role in our method are formulated in  \S4 and proved in
\S5.
 \S2 is devoted to    precise setting of the problem. In \S3 we recall necessary homogenization attributes, among them, first of all, 
 the so called cell problems and  the homogenized matrix which is defined with the help of the solutions to the cell problem.  In locally periodic setting the cell problem depends on the slow variable $x$ as a parameter. So  one should take into account some properties of the cell problem solutions with respect to the variable $x$. We recall these properties, known from the previous papers on locally periodic homogenization, also in \S3.

\textbf{1.3. About method.}
The present paper can be viewed as following in the footsteps of  \cite{Zh1} 
 in that it relies upon the so-called  "modified method of the first approximation" with the usage of  shift. The same is true for the  above-mentioned  paper \cite{PT}. On the contrary, in \cite{BS} and also in   \cite{Zh05} and \cite{BS05}, the authors applied the spectral approach that seems to be 
 tightly coupled with the assumption of periodicity in coefficients, because one of the components of this approach is the Floquet--Bloch transformation   which serves well only in periodic setting. The modified method of the first approximation was proposed by V.V.Zhikov \cite{Zh1}
 as an alternative approach to prove operator-type homogenization estimates of the type (\ref{0.3}) and (\ref{0.6}); it 
 turned to be universal in different setups with periodic, locally periodic, quasiperiodic or multiscale coefficients.
The method has developed since 2005  in applications to various problems (we refer, e.g.,  to \cite{ZhP05}--\cite{PT17} and, in particular, to the overview \cite{UMN} where other  references are given); there have appeared two versions of it: 
one with the usage of the pure shift,  and another with the usage of the Steklov smoothing operator applying the shift implicitly.

 Recently \cite{P20}, by the modified method of the first approximation, 
 the  estimate (\ref{0.5}) was  proved for the problem  (\ref{0.1}) in periodic setting, in so doing, the version with Steklov's smoothing was chosen. 
 Now we address the locally periodic setting and prove the estimate of the type (\ref{0.5}), choosing the original version of this method (coming  straight from \cite{Zh1}) with the pure shift. Steklov's smoothing arises here only at the  ultimate step of the proof and participates in the final formulation of the result.
 
 For the sake of simplicity we restrict ourselves to the scalar case.
 We consider the classical
diffusion equation of the type (\ref{0.1}), but with a locally periodic  not necessarily symmetric  diffusion matrix, given in the whole space. The
obtained result admits generalizations.
 Although we deal with the classical diffusion equation, the maximum principle
or its corollaries, valid in the scalar case, are not used in our constructions, and so
 the result also carries over to vector models, including, e.g., the elasticity theory system.
 
 It is also worth noting that, once the  estimate (\ref{0.6}) (or 
   mentioned above its more general counterpart with Steklov's smoothing in the corrector) in the operator $(L^2\to H^1)$-norm with order $\e$ remainder is verified, 
 the estimate of the type (\ref{0.5})
 in the operator $(L^2\to L^2)$-norm with order $\e^2$ remainder is surely guaranteed by the method we demonstrate here.
 
 Our addressing 
 the estimates of the type 
 (\ref{0.5}) for locally periodic elliptic operators 
 appears as a response to publications \cite{Se0}, \cite{Se}. 
    We are aimed  to show here that the shift method proposed by Zhikov in \cite{Zh1} is quite 
effective towards   
     this issue either. 
   That confirms once more  the high potential of the shift  method.
   
   As  a by-product of this paper we obtain an alternative proof for  the results of \cite{P20}, \cite{P20A} where resolvent approximations   of elliptic selfadjoint or nonselfadjoint differential
 operators with  periodic coefficients  were studied.
The proof of the present paper has an advantage over those given in \cite{P20}, \cite{P20A} because it can be extended,  in a natural fashion  (see, e.g., \cite{PT}, \cite{PaT15} or \cite{UMN}), to the problems with multiscale coefficients where reiterated homogenization takes place.

\section{Problem setup}


Consider the following elliptic equation in the whole space
 $\RR^d$:
 \beq\label{eq1.1}
\ds{u^\ve\in H^1(\RR^d),\quad A_\ve u^\ve+ u^\ve=f,\quad
f\in L^2(\RR^d),} \atop\ds{A_\ve =-\div\, a^\ve(x)\nab},\eeq
where the matrix
$a^\ve(x)$ is locally periodic. The latter means that
\beq\label{eq1.2}
a^\ve(x)=a(x,\e^{-1}x
),\quad \e>0,
 \eeq
and $a(x,\cdot)$ is  $Y$-periodic, where the periodicity cell is the unit cube $Y=[-1/2,1/2)^d$. Besides, the matrix function $a(x,y)$ satisfies the Carat\'{e}odory 
condition (with continuity in $x$ and measurability in $y$) which garantees  measurability of the locally periodic function
$a^\ve(x)$. According to (\ref{eq1.2}), $a^\ve(x)$ is  rapidly oscillating,  as the small parameter $\ve$ goes to zero, but this oscillation clearly is not periodic.  

Moreover, we require the following conditions  on $a(x,y)$:
\beq\label{eq1.3}
|a(x,y)-a(x',y)|\le c_L|x'-x|
\eeq
for all $x, x'\in \RR^d$ and a.e. $y\in  Y$, that is the Lipschitz continuity of  $a(x,y)$ with respect to the first variable;
\beq\label{eq1.4}
\lambda |\xi|^2\le {a}\xi\cdot\xi
,\quad {a}\xi\cdot\eta
\le\lambda^{-1}|\xi|\,|\eta|
\quad\forall \xi,\eta\in\R^d
\eeq
for some $\lambda>0$. The matrix $a$ is not necessarily symmetric.

Equation (\ref{eq1.1}) is related to the homogenized equation
\beq\label{eq1.5}
u\in H^1(\RR^d),\ A_0u +u=-\div\, a^0(x)\nabla u+u=f,
\eeq
where the matrix $a^0(x)$ depends only on the  \ "slow"\ variable $x$, satisfies the conditions of the type (\ref{eq1.3}), (\ref{eq1.4}) and may be found through the known procedure by solving auxiliary problems on the periodicity cell $Y$ (see (\ref{eq2.1}), (\ref{eq2.2})).

Equations (\ref{eq1.1}) and (\ref{eq1.5}) are understood in the sense of the theory of distributions on $\RR^d$; they are uniquely solvable  for any  right-hand side 
 function $f\in L^2(\RR^d) $ (even for
any function $f$ from a wider space $H^{-1}(\RR^d)$ that is a dual space to $H^1(\RR^d)$), and the uniform  (in $\e$) energy estimate
\beq\label{ener}
\|u^\ve\|_{ H^1(\RR^d)}\le c\|f\|_{ L^2(\RR^d)},\quad c=const(d),
\eeq
is fulfilled.
 Thus, one can speak about the resolvents $(A_\ve+1)^{-1} $ and $(A_0+1)^{-1}$ acting in $L^2(\RR^d)$.
 The question arises in what sense these resolvents are close to each other.
The long-known result says  that, for any right-hand side 
 function $f\in L^2(\RR^d) $, the solutions of the original and the homogenized problems are connected by the
 strong convergence 
$u^\ve\to u$ in $L^2(\RR^d)$. 
In  the operator terms, this means the strong resolvent convergence
\[
(A_\ve+1)^{-1}\to (A_0+1)^{-1}\quad \mbox{in}\quad L^2(\RR^d).
\]
The latter can be strengthened to the uniform resolvent convergence 
with the following sharp with respect to the order
estimate for the convergence rate
\beq\label{eq1.7}
\|(A_\ve+1)^{-1}- (A_0+1)^{-1}\|_{L^2(\RR^d)\to L^2(\RR^d)}\le C\ve,
\eeq
where the constant $C$ depends only on the dimension $d$ and  the constants of the conditions (\ref{eq1.3}), (\ref{eq1.4}). So, the zero approximation for the resolvent  $(A_\ve+1)^{-1}$ of the original operator in the operator $L^2$-norm is  the resolvent  $(A_0+1)^{-1}$ of the homogenized operator.

 If the resolvent $(A_\ve+1)^{-1}$ is regarded as an operator from 
 $L^2(\RR^d)$ to $H^1(\RR^d)$, then for its approximations we take the sum of the constructed  zero approximation and some correcting operator, i.e., $(A_0+1)^{-1}+\e\K_\e$, where $\K_\e$ is defined in (\ref{eq3.5o})
 and (\ref{eq3.5c}).
  Then
\beq\label{eq1.8}
\|(A_\ve+1)^{-1}- (A_0+1)^{-1}-\e\mathcal{K}_\ve\|_{L^2(\RR^d)\to H^1(\RR^d)}\le C\ve
\eeq
with the constant $C$  of the same type, as in  (\ref{eq1.7}).

 The operator estimates (\ref{eq1.7}) and (\ref{eq1.8}) were proved in
 \cite{PT} (see also \cite{UMN}). To obtain 
 the correcting operator $\K_\e$ involved in (\ref{eq1.8}), a  perturbated family of 
 operators $A_\e^\o$ with a shifted matrix
 $a^\e_\o$ was introduced, $\o$ being a shift parameter. 
 Then 
 the auxiliary, averaged over the shift parameter $\o$, 
  $H^1$-estimate for the difference between $(A_\e^\o+1)^{-1}f$ and
  the appropriate approximation
was established. From this the  approximations for  $(A_\ve+1)^{-1}f$
 naturally arose, in which  Steklov`s smoothing was comprised, with the
 estimate (\ref{eq1.8}) following as a corollary. In \S5 we reproduce in details this derivation of (\ref{eq1.8}) relying on the idea of shift, because the elements of this proof are systematically   used further when we address the $L^2$-approximations of  $(A_\e+1)^{-1}$ with $\e^2$ remainder estimates and 
 seek an appropriate corrector $\mathcal{C}_\e$ such that
\beq\label{eq1.9}
\|(A_\ve+1)^{-1}- (A_0+1)^{-1}-\e\mathcal{C}_\ve\|_{L^2(\RR^d)\to L^2(\RR^d)}\le C\ve^2,
\eeq
with the constant $C$  of the same type, as in  (\ref{eq1.7}). The corrector $\mathcal{C}_\ve$ is constructed in \S5, it turns to be a sum of several terms of different structure, namely,
 \beq\label{eq1.10}
 \mathcal{C}_\ve=\K_\e +\tilde{\K}_\e^*-\mathcal{L}-\mathcal{M}_\e.
 \eeq
The operators  $\K_\e $ and $\tilde{\K}_\e$ contain oscillating locally periodic factors dependant on $\e$ (similarly, as $\K_\e $ in periodic setting, see (\ref{0.5}) and (\ref{0.5c})). The operator $\mathcal{L}$ does not depend on $\e$ at all which is indicated in notation. As for $\mathcal{M}_\e$, there is no oscillating factors in its structure and  $\e$ participates in it as a parameter of 
smoothing involved in the operator. 
 
 It is  worth noting that to construct all the four terms in $\mathcal{C}_\ve$ 
one  uses only the resolvent $(A_0+1)^{-1}$ of the homogenized operator and 
the solutions of  the two cell problems (\ref{eq2.1}) and (\ref{cps}) given below (and of only one cell problem (\ref{eq2.1}) if the matrix $a$ is symmetric),
no other auxiliary problems are needed to this end.
 In the case of the periodic selfadjoint setting, the operator $ \mathcal{C}_\ve$  written in (\ref{eq1.10}) reduces to the corrector  in (\ref{0.5}) consisting of only two terms with $\tilde{\K}_\e=\K_\e$ and  zero $\mathcal{L}$, $\mathcal{M}_\e$.

The precise formulation of the 
main result, that is    the $L^2$-estimate (\ref{eq1.9})  with the correcting term (\ref{eq1.10}), 
is given in Theorem \ref{Th5.1}.

\section{Homogenization attributes}

\textbf{3.1. Cell problems.}
Consider  periodic problems on the unit cube $Y=[-\f{1}{2},\f{1}{2})^d$
\beq\label{eq2.1}
\ds{
N^j(x,\cdot)\in H^1_{per}(Y),\ \div_y [a(x,y)(e^j+\nabla_y N^j(x,y))]=0,\ \langle
N^j(x,\cdot)\rangle=0,}\atop\ds{
j=1,\ldots,d,}
\eeq
where the variable $x$ plays the role of a parameter. We use here the notation:
$e^1, \dots , e^d$ is a canonical basis in $\R^d$,
$H_{\rm per}^1(Y)$ is the Sobolev space of 1-periodic functions,
 $$\langle \cdot\rangle=\langle \cdot\rangle_Y=\int\limits_Y\,\cdot\, dy.$$ 
 The equation (\ref{eq2.1})  can be understood either in the sense of
 distributions on $\RR^d$ or in the sense of the integral identity on the periodicity cell. The latter means 
 \[
 \langle a(x,\cdot)(e^j+\nabla_y N^j(x,\cdot))\cdot\nab \varphi\rangle=0\quad \forall \varphi\in C^\infty_\per(Y).
 \]
 These two formulations  
 of the cell problem are equivalent which will be taken into account later. 
 
The homogenized matrix $a^0$ is defined 
 in terms of the solutions to the cell problems by the following relations:
\beq\label{eq2.2}
a^0(x)e^j=\langle a(x,\cdot)(e^j+\nabla_y N^j(x,\cdot))\rangle,\, j=1,\ldots,d,
\eeq
and depends on the "slow"\ variable $x$.
The ) 
  matrix $a^0(x)$ inherits the Lipschitz continuity in $x$ from the original matrix $a(x,y)$ (see below Lemma \ref{locper}),
  thereby, the ellipticity theory 
 yields the estimate
 \beq\label{ells}
\|u \|_{ H^2(\R^d)}\le c \|f\|_{ L^2(\R^d)}, \quad 
c=cost(\lambda), 
\eeq
for the solution to (\ref{eq1.5}).

We list now  some properties of the solutions to the cell problem.

\begin{lemma}\label{locper}

Let $N^j(x,y)$ be the solutions of the problems (\ref{eq2.1}), and let $a^0(x)$ be the matrix defined in (\ref{eq2.2}). Then:

i) $\|N^j(x,\cdot)\|_{H^1_{per}(Y)}\le c$ for all $x$,
 where $c=const(d,\lambda)$;

ii) $N^j(x,y)$ is a Lipschitz continuous function with respect to $x$ with values in  $H^1_{per}(Y)$, and its  Lipschitz constant depends only on the constants $c_L$, $\lambda$ from (\ref{eq1.3}), (\ref{eq1.4});

iii) the matrix $a^0(x)$ is  Lipschitz continuous with the Lipschitz constant depending only on the constants $c_L$, $\lambda$ from (\ref{eq1.3}), (\ref{eq1.4});

 iv) there exists the gradient $\nabla_x N^j$ such that $\|\nabla_x N^j(x,\cdot)\|_{H^1_{per}(Y)^d}\le c$ with the constant $c=const(d,\lambda,c_L)$.
\end{lemma}
\textbf{Proof} of the properties  $ii),iii)$  is given, e.g., in \cite{PT}, \cite{PAsAn}, \cite{PaT15}.
The property $iv)$ is a corollary of 
$i), ii)$ due to Rademacher's theorem: every Lipschitz function belongs to $W^{1,\infty}$-space. The property $i)$ follows from the energy estimate for the solution of the cell problem and our assumptions on the matrix $a(x,y)$.


\textbf{3.2. Adjoint problems.}
 Let $A_\e^*$ be the adjoint of
  $A_\e$ and consider the problem
 \begin{equation}\label{1s}\ds
u_\ve \in H^1 (\R^d), \quad 
{A_{\ve}^*v^\ve+v^\ve=h,\quad h{\in} L^2(\R^d),}
\atop\ds{A_{\ve}^*=-\div\,(a^\ve(x))^*\nab,\quad (a^\ve(x))^*=a^*(x,\ve^{-1}x),}
\end{equation}
where 
$a^*(x,y)$ is the transposed matrix to $a(x,y)$.

It is known that the homogenized equation for
(\ref{1s}) 
 will be  
\beq\label{homs}
v \in H^1 (\R^d), \quad 
 A^*_0 v+v=-\div \,{(a^0)^*}\nab v+v=h,
\eeq
where $A^*_0$ is the adjoint of $A_0$ and has the matrix   $(a^0)^*$ transposed to ${a^0}$.
Thus,
\beq\label{hom}
(a^*)^0 =(a^0)^*.
\eeq

 In the case of the adjoint equation, the counterpart of  the cell problem  (\ref{eq2.1}) will be
\beq\label{cps}\ds{
\tilde{N}^j(x,\cdot)\in H_{\rm per}^1(\Box),\quad\div_y a^*(x,y)(e^j+\nab_y \tilde{N}^j)=0, }\atop\ds
{\langle \tilde{N}^j(x,\cdot)\rangle=0,\quad 
j=1,  ..., d.}
\eeq
Its solutions  generate formally  the homogenized matrix for the equation (\ref{1s}) 
 through the formula similar to (\ref{eq2.2}), and so $\tilde{N} ^j$  are connected with the matrix ${a^0}$: 
\beq\label{13s}
{(a^0)^*e^j}=\langle a^*(e^j+\nab \tilde{N}^j)\rangle,\quad j=1, \dots ,d,
\eeq
where (\ref{hom}) is taken into account.

\textbf{3.3. Shifting and smoothing operators.} Given $\e\in(0,1)$ and $\o\in Y=[-1/2,1/2)^d$, we use the notation
\[
S^\e_\o\varphi(x)=\varphi(x+\e\omega),
\]
\beq\label{st}
S^\e\varphi(x)=\int\limits_Y\varphi(x-\e\omega)\,d\omega 
\eeq
for the shift operator and the Steklov average, the latter is also referred to as the Steklov smoothing operator.
In our method these operators play the key role. We list 
 here their properties that will be used in the sequel:
\beq\label{sh1}
\|S^\e_\o\varphi-\varphi\| \le c\e\|\nab \varphi\|,
\eeq
\beq\label{st1}
\|S^\e\varphi-\varphi\| \le c\e\|\nab \varphi\|,
\eeq
\beq\label{sh2}
\|S^\e_\o\varphi-\varphi\|_{H^{-1}(\rd)} \le c\e\|\varphi\|,
\eeq
where $\|\cdot\|$ denotes the norm $\|\cdot\|_{L^2(\rd)}$ and the constant in the right-hand side depends only on the dimension $d$.
The estimate (\ref{st1}) may be sharper if the function $\varphi$ is more regular, more precisely, if its second gradient $\nab^2 \varphi$ belongs to $L^2(\rd)$:
\beq\label{st2}
\|S^\e\varphi-\varphi\| \le c\e^2\|\nab^2 \varphi\|, \quad c=const(d).
\eeq

The property (\ref{sh2}) is a  corollary of (\ref{sh1}), by duality arguments. The other properties follow from the Taylor
 formula  with the remainder term in its simplest forms applied to write the difference $\varphi(x+\e\omega)-\varphi(x)$ in terms of derivatives.
\section{Shifted first approximation}

\textbf{4.1. A perturbated family of problems.} In classical homogenization the expression 
\beq\label{eq3.0}
u(x)+\ve U(x,y):=u (x)+
\ve N^j(x,y)\f{\pa u(x)}{\pa x_j}=
u (x)+\e N(x,y)\cdot \nab u(x),\quad y=\f{x}{\ve},
\eeq
is commonly called the first approximation of the solution $u^\e(x)$ to (\ref{eq1.1}), in its turn, the term $\ve N^j(x,\f{x}{\ve})\f{\pa u(x)}{\pa x_j}$ is called a corrector. The solutions of the cell problem (\ref{eq2.1}) and the homogenized eqution (\ref{eq1.4}) are involved here, they are not too much regular under our assumptions,  so the corrector 
$\ve N^j(x,\f{x}{\ve})\f{\pa u(x)}{\pa x_j}$ is not necessarily from $ H^1(\RR^d)$. Thereby, 
the function (\ref{eq3.0}), generally, does not belong to $H^1(\RR^d)$,  thus,  it cannot approximate the solution $u^\e(x)$ in the norm of
$H^1(\RR^d)$. In what follows, we show how to overcome this difficulty by 
 modifying the 
concept of the first approximation.

We consider  a family of  perturbated  problems 
\beq\label{eq3.1}
\ds{u^\ve_\omega \in H^1(\RR^d),\quad A^\omega_\ve u_\omega^\ve+ u^\ve_\omega =f,\quad
f\in L^2(\RR^d),} \atop\ds{A^\omega_\ve =-\div a^\ve_\omega (x)\nab},
\eeq
with the shifted matrix
$$
a^\e_\omega(x):=a\left(x,\e^{-1}x
+\omega\right),\ \omega\in \RR^d.
$$
Clearly, taking $\omega=0$ in (\ref{eq3.1}), we come to the original equation (\ref{eq1.1}).

It is evident that the solution of the periodic problem (\ref{eq2.1}), with
the shifted matrix $a(x,y+\omega)$ instead of $a(x,y)$, is obtained from $N^j(x,\cdot)$ by shifting in the argument by
$\omega$, that is
$N^j(x,\cdot+\omega)$. However, the calculation of the homogenized matrix for 
(\ref{eq3.1}) according to  (\ref{eq2.2}) gives the same $a^0(x)$, just like for  (\ref{eq1.1}); 
therefore, (\ref{eq3.1}) relates to the homogenized equation defined in (\ref{eq1.4}).
 Then
\beq\label{eq3.2}
u(x)+\ve U(x,y+\omega):=u(x)+\ve N^j(x,y+\omega)\f{\pa u(x)}{\pa x_j},\ y=\frac{x}{\ve},
\eeq
turns to be the first approximation for the solution of (\ref{eq3.1}).
We have a shifted corrector in (\ref{eq3.2}). As a function of two variables $x$ and $\omega$, $U(x,\e^{-1}x+\omega)$ belongs to 
$L^2(\RR^d\times Y)$; the same is valid for its gradient in $x$.
Moreover, the following estimate holds:
\beq\label{estcor}
\int\limits_{Y\times\RR^d}
 \l(|U(x,\frac{x}{\ve}+\omega)|^2+|\e\nab_x U(x,\frac{x}{\ve}+\omega)|^2
 \r)dx\,d\o\le c \|f\|^2_{L^2(\rd)}, \quad  c=const(d,\lambda,c_L).
\eeq
Indeed,  consider, for example,  the gradient 
$$\ds{
\e\nab_x \l( U(x,\frac{x}{\ve}+\omega)\r)=\nabla_y N^j(x,y+\omega)\f{\pa u(x)}{\pa x_j}+ \e\nabla_x N^j (x,y+\omega)\f{\pa u(x)}{\pa x_j}+}
\atop\ds{\ve N^j(x,y+\omega)\nabla\f{\pa u(x)}{\pa x_j},\ y=\frac{x}{\ve}.}
$$
We have
$$
\int\limits_{Y\times\RR^d}\l|\e\nab_x  U(x,\frac{x}{\ve}+\omega)
\r|
^2dxd\omega\le
$$
$$
\sum_k\int\limits_{Y\times\RR^d}|b_k(x,\frac x\ve+\omega)|^2|\Phi (x)|^2dxd\omega=
\sum_k\int\limits_{\RR^d}
\l(\int\limits_{Y}|b_k(x,\frac x\ve+\omega)|^2d\omega\r)|\Phi (x)|^2
 dx\le
 $$
 $$\int\limits_{\RR^d}|\Phi (x)|^2
 dx\sum_k \sup_x\|b_k(x,\cdot)
\|_{L^2_\per(Y)}
 \le C\int\limits_{\RR^d}|\Phi (x)|^2
 dx.
$$
Here
$\Phi(x)=|\nabla u (x)|+|\nabla^2 u (x)|$ and  
$\Phi\in L^2(\RR^d)$  by the estimate (\ref{ells});
 the terms $b_k(x,y)$ are formed  of 
$$ N^j(x,y),\,
\nabla_x N^j(x,y),\, \nabla_y N^j(x,y),$$
for which, by Lemma \ref{locper}, 
 the boundedness property
$\sup_x\|b_k(x,\cdot)
\|_{L^2_\per(Y)}\le C$ is guaranteed.
 Thus, the estimate (\ref{estcor}) 
is verified.

We have just actually applied and proved 
the following 
 \begin{lemma}\label{LemShift} Suppose $\varphi{\in}L^1(\R^d)$, $b{\in}L^\infty( \RR^d,L^1_\per(\Box))$, and define $b_\e^\omega(x)=b(x,\e^{-1}x+\omega)$. Then, as a function of
 two variables $x$ and $\omega$,
  the product
  $\varphi(x)\, b_\e^\omega(x)$  belons to the space $ L^1(\R^d\times Y)$ and the following estimate holds
 \begin{equation}\label{eq3.4}
\|b_\e^\omega\varphi\|_{L^1(\rd\times Y)}\le\|b_\e^\omega
\|_{L^\infty( \RR^d,L^1_\per(\Box))}\|\varphi\|_{L^1(\rd)}.
\end{equation}
\end{lemma}

The above arguments show 
 that the approximation (\ref{eq3.2})
and its gradient with respect to $x$
belong to 
$L^2(\RR^d\times Y)$ as the functions of
 two variables $x$ and $\omega$; besides, 
their $L^2$-norms on the product $\RR^d\times Y$ are uiniformly in $\e$ bounded.
Moreover, the function (\ref{eq3.2})
 approximates the solution to (\ref{eq3.1}) in
the following $\omega$-averaged sense.

\begin{lemma}\label{LemAv}
Let $u^\ve_\omega(x)$  be the solution to (\ref{eq3.1}) and let $u(x)+\ve U(x,\e^{-1}x+\omega)$ be the corresponding first approximation defined in (\ref{eq3.2}). Then 
the $\omega$-averaged inequality
\beq\label{eq3.5}
\int\limits_{Y} \|u^\ve_\omega (\cdot)-u(\cdot)-\ve U(\cdot,\e^{-1}\cdot+\omega)\|^2_{H^1(\RR^d)}d\omega\le C\ve^2 \|f\|^2_{L^2(\RR^d)}
\eeq
holds, where   the constant $C$ depends only on the dimension $d$ and the constants $c_L$, $\lambda$ from (\ref{eq1.3}), (\ref{eq1.4}).
\end{lemma}

There is another version of the $\omega$-averaged $H^1$-estimate which show 
 a direct relationship between the shifted first approximation defined in (\ref{eq3.2}) and the solution of the original problem (\ref{eq2.1}).
 
 \begin{lemma}\label{LemAv1}
Let $u^\ve(x)$  be the solution to (\ref{eq2.1}) and  let $u(x)+\ve U(x,\e^{-1}x+\omega)$ be the shifted first approximation defined in (\ref{eq3.2}). Then 
the $\omega$-averaged inequality
\beq\label{eq3.5a}
\int\limits_{Y} \|u^\ve (\cdot+\e\omega)-u(\cdot)-\ve U(\cdot,\e^{-1}\cdot+\omega)\|^2_{H^1(\RR^d)}d\omega\le C\ve^2 \|f\|^2_{L^2(\RR^d)},
\eeq
holds, where   the constant $C$ is of the same type, as in (\ref{eq1.8}).
\end{lemma}

As a corallary of Lemma \ref{LemAv1}, we obtain
 \begin{lemma}\label{LemAv2}
Let $u^\ve(x)$  be the solution to (\ref{eq2.1}) and let 
 $u(x)+\ve U(x,\e^{-1}x)=u(x)+\ve N(x,\e^{-1}x)\cdot \nab u(x)$
 be the  first approximation defined in (\ref{eq3.0}). Define 
\beq\label{eq3.5c}
K_\e(x)= \int_Y N(x-\e\o,\e^{-1}x)\cdot \nab u(x-\e\o)\,d\,\o
\eeq
that is a smoothed corrector. Then the following $H^1$-estimate 
\beq\label{eq3.5b}
 \|u^\ve -u-\ve K_\e\|_{H^1(\RR^d)}\le C\ve \|f\|_{L^2(\RR^d)}
\eeq
holds, where   the constant $C$ is of the same type, as in (\ref{eq3.5}).
\end{lemma}

These lemmas are proved in the next section. 

Define the operator  $\K_\e:\ld \to H^1(\RR^d)$ as follows
\beq\label{eq3.5o}
 \K_\e f:=K_\e
\eeq
with $K_\e$ given in (\ref{eq3.5c}). Then (\ref{eq3.5b}) is equivalent to (\ref{eq1.8}).

\section{ Proof of  $\omega$-averaged estimates}

1$^\circ$ Let $v^\e(x)$ denote the approximation (\ref{eq3.0}). We start with an analysis
of this approximation and its discrepancy in the original equation (\ref{eq1.1}). 

Calculating the gradient $\nab v^\e(x)$, 
we compare the fluxes $a^\e(x)\nab v^\e(x)$ and $a^0(x)\nab u(x)$:
\beq\label{grad}
\nabla v^\ve (x){=}\nabla u(x) +\nabla_y N^j(x,y)\f{\pa u(x)}{\pa x_j}{+}\ve \nabla_x N^j (x,y)\f{\pa u(x)}{\pa x_j}{+}\ve N^j(x,y)\nabla  \f{\pa u(x)}{\pa x_j},\ y=\frac{x}{\ve},
\eeq
\beq\label{eq3.6}\ds{
a(x,y)\nabla v^\ve (x)- a^0 (x)\nabla u (x)=a(x,y)(e^j+\nabla_y N^j(x,y))\f{\pa u(x)}{\pa x_j}
}\atop\ds{
- \langle a(x,\cdot)(e^j+\nabla_y N^j(x,\cdot))\rangle_Y\f{\pa u(x)}{\pa x_j}+r_{1,\ve},\ y=\frac{x}{\ve},}
\eeq
with
\beq\label{eq3.7}
r_{1,\ve} = \ve a(x,y)\nabla_x N^j(x,y)\f{\pa u(x)}{\pa x_j}+\ve  a(x,y)N^j(x,y)\nabla\f{\pa u(x)}{\pa x_j},\ y=\frac{x}{\ve},
\eeq
where we have used the definition of the homogenized matrix $a^0 (x)$ (see
(\ref{eq2.2})). 
 Hence
\beq\label{eq3.8}
a(x,y)\nabla v^\ve(x)-a^0(x)\nabla u (x)=g^j(x,y)\f{\pa u(x)}{\pa x_j}+ r_{1,\ve},\ y=\frac{x}{\ve},
\eeq
with
\beq\label{eq3.80}
g^j(x,y)= a(x,y)(e^j+\nabla_y N^j(x,y))-\langle a(x,\cdot)(e^j+\nabla_y N^j(x,\cdot))\rangle_Y.\eeq
Evidently, the periodic (over $y$) vector
$g^j(x,y)$ is such that 
\beq\label{solg}
\Div_y g^j(x,y)=0,\quad \langle g^j(x,\cdot)\rangle_Y=0.
\eeq
Besides, by Lipschitz continuity properties of the matrix $a(x,y)$ and  the  solutions $N^j(x,y)$ (see (\ref{eq1.3}) and Lemma  \ref{locper}), the function
$g^j(x,y)$
turns  to be Lipschitz continuous over $x$ with values in $L^2_{per}(Y)$.
This property of $g^j(x,y)$ allows us to prove (see the proof in \cite{PT}, \cite{PaSt}, \cite{PAsAn}, \cite{PaT15}).

\begin{lemma}\label{locper1}
There exists a skew-symmetric matrix $G^j(x,y)$ ($G^j_{ik}=-G^j_{ki}$)
such that
\beq\label{eq3.9}
g^j(x,y)=\Div_y G^j(x,y)
\eeq 
and 
$$
\|G^j(x,\cdot)\|_{H^1_{per}(Y)^{d\times d}}\le c\|g^j(x,\cdot)\|_{L^2_{per}(Y)^d},\quad c=c(\lambda,c_L).
$$
Moreover, $G^j(x,\cdot)$ is Lipschitz continuous in $x$ with values in $H^1_{per}(Y)^{d\times d}$, thereby,
its gradient in $x$ exists for a.e. $x\in \RR^d $, and the estimate
$$
\|\nabla_x G^j_{ik}(x,\cdot)\|_{H^1_{per}(Y)^d}\le c,\quad c=c(\lambda,c_L).
$$
holds for  a.e. $x\in \rd$.
\end{lemma}

 Note that the equality (\ref{eq3.9})  as well as (\ref{solg})$_1$ can be understood in two 
 ways, the same as the equation (\ref{eq2.1}).

In view of (\ref{eq3.9}), we  write 
\beq\label{eq3.10}\ds{
g^j(x,y)\frac{\partial u(x)}{\partial x_j}=\ve\Div \left(G^j(x,y)\frac{\partial u(x)}{\partial x_j}\right)
}\atop\ds{
-\ve \left(\Div_x G^j(x,y)\frac{\partial u(x)}{\partial x_j}+G^j(x,y)\nabla\frac{\partial u(x)}{\partial x_j}\right),\ y=\frac{x}{\ve},}
\eeq
where the first term in the right-hand side is a solenoidal vector. Indeed, 
\beq\label{solG}
\int\limits_{\RR^d}\Div \left(G^j (x,\frac{x}{\ve})\frac{\partial u(x)}{\partial x_j}\right)\cdot\nabla\varphi (x) dx=-\int\limits_{\RR^d}G^j (x,\frac{x}{\ve})\cdot\nabla^2\varphi (x)\frac{\partial u(x)}{\partial x_j} dx=0
\eeq
for any $\varphi\in C^\infty_0(\RR^d)$,
because the matrices $\nabla^2\varphi$ and $G^j$ are symmetric and skew-symmetric respectively, thereby,
$G^j\cdot \nabla^2\varphi=0$ pointwise.

Thus, from (\ref{eq3.6})--(\ref{eq3.8}) and (\ref{eq3.10}) we derive
\beq\label{eq3.11}
\ds{
\Div[a(x,y)\nabla v^\ve(x)-a^0(x)\nabla u(x)]
=\Div\bigg(\ve a(x,y)\nabla_x N^j(x,y)\frac{\partial u(x)}{\partial x_j}}
\atop\ds{+\ve a(x,y) N^j(x,y)\nabla \frac{\partial u(x)}{\partial x_j}
-
\ve \Div_x G^j(x,y)\frac{\partial u(x)}{\partial x_j}-\ve G^j(x,y)\nabla\frac{\partial u(x)}{\partial x_j}\bigg),\ y=\frac{x}{\ve}.}
\eeq

Since
\beq\label{eq3.120}\ds{A_\ve (v^\ve -u^\ve)+(v^\ve - u^\ve)=A_\ve v^\ve+v^\ve-A_0 u - u=}
\atop\ds{
=-\Div[a(x,\frac{x}{\ve})\nabla v^\ve(x)-a^0(x)\nabla u(x)]+v^\ve - u,}
\eeq
we obtain the equation
\beq\label{eq3.12}
A_\ve w^\ve + w^\ve=f_\ve+\Div F_\ve
\eeq
for $$w^\ve=v^\ve -u^\ve,
$$
where (see  (\ref{eq3.120})) and (\ref{eq3.11}))
$$
f_\ve=\ve N^j(x,\frac{x}{\ve})\f{\pa u(x)}{\pa x_j},
$$
$$
F_\ve=
-\ve\bigg( a(x,y)\nabla_x N^j(x,y)\frac{\partial u(x)}{\partial x_j}+  a(x,y) N^j(x,y)\nabla\frac{\partial u(x)}{\partial x_j}
$$
$$
- \Div_x G^j(x,y)\frac{\partial u(x)}{\partial x_j}-G^j(x,y)\nabla\frac{\partial u(x)}{\partial x_j}\bigg),\ y=\frac{x}{\ve}.
$$
The energy inequality
$$
\|w^\ve\|_{H^1(\RR^d)}\le c_0 (\|f_\ve\|_{L^2(\RR^d)}+\|F_\ve\|_{L^2(\RR^d)}),\ c_0=const(\lambda),
$$
is fulfilled for the equation (\ref{eq3.12}).
Taking into account the expressions for the right-hand side functions
$f_\e$ and $F_\e$ in (\ref{eq3.12}), we can write
\beq\label{eq3.13}
\|v^\ve - u^\ve\|^2_{H^1(\RR^d)}\le c_0\ve^2\sum_k\int\limits_{\RR^d}|b_k(x,y)|^2|\Phi (x)|^2dx,\quad y=\f{x}{\ve}.
\eeq
Here $\Phi(x)=|\nabla u (x)|+|\nabla^2 u (x)|$, and the terms $b_k(x,y)$ are formed of the functions
\beq\label{eq3.14}
G^j(x,y), \,\nabla_x G^j(x,y),\,N^j(x,y),\,\nabla_x N^j(x,y).
\eeq

2$^\circ$
Let us try to exclude the factors $b_k(x,y)$ from the integrals in the estimate (\ref{eq3.13}).
To this end, we address the perturbated problem (\ref{eq3.1}), for which the first approximation is defined in (\ref{eq3.2}), and write the counterpart of the estimate (\ref{eq3.13}) relating to (\ref{eq3.1}). 
Let, for brevity, $v^\ve_\omega(x)$ denote the first approximation  defined in (\ref{eq3.2}). Then, according to (\ref{eq3.13}), we can write
\[
\|u^\ve_\omega-v^\ve_\omega \|^2_{H^1(\RR^d)}\le c_0 \ve^2\sum_k\int\limits_{\RR^d} \left|b_k\left(x,\frac x\ve+\omega\right)\right|^2|\Phi(x)|^2dx,
\]
 next integrate over $\omega\in Y$  
$$
\int\limits_{Y}\|u^\ve_\omega(\cdot)-v^\ve_\omega(\cdot)\|^2_{H^1(\RR^d)}d\omega=
$$
$$
=\int\limits_{Y}\int\limits_{\RR^d}(|u^\ve_\omega(x)-v^\ve_\omega (x)|^2+|\nabla u^\ve_\omega (x)-\nabla v^\ve_\omega (x)|^2)dxd\omega\le
$$
$$
\le c_0\ve^2\sum_k\int\limits_{Y}\int\limits_{\RR^d}|b_k(x,\frac x\ve+\omega)|^2|\Phi (x)|^2dxd\omega
\le c_0\ve^2\int\limits_{\RR^d}|\Phi (x)|^2dx\cdot\sum_k\sup_x\int\limits_{Y}|b_k(x,\frac x\ve+\omega)|^2d\omega,
$$
and finally deduce
\beq\label{eq3.16}
\int\limits_{Y}\|u^\ve_\omega(\cdot)-v^\ve_\omega(\cdot)\|^2_{H^1(\RR^d)}d\omega
\le C\ve^2 
\|f\|^2_{L^2(\RR^d)},\ C=const (d,\lambda,c_L).
\eeq
Here at the last steps we have applied Lemma \ref{LemShift}, the elliptic estimate  (\ref{ells}), and also the estimate
$$
\int\limits_{Y}|b_k(x,\omega)|^2 d\omega\le c,\, c=const (d,\lambda,c_L),
$$
 which is valid due to the properties of the functions (\ref{eq3.14}) established earlier (see Lemma \ref{locper} and Lemma  \ref{locper1}).

The estimate (\ref{eq3.16}) is equivalent to the desired estimate (\ref{eq3.5}): it is enough only to explicate the notation $v^\e_\o$. Lemma \ref{LemAv} is proved.

\bigskip
3$^\circ$ We proceed now to Lemma \ref{LemAv1} and compare 
the solution to the problem (\ref{eq3.1}) and the function
$u^\ve (x+\ve\omega)$ which is the solution of the equation
\beq\label{eq3.18}
-\Div \l(a(x+\ve\omega,\e^{-1}x
+\omega)\nabla u^\ve (x+\ve\omega)\r)+ u^\ve (x+\ve\omega)=f(x+\ve\omega).
\eeq
Setting
$w^\ve_\omega (x)=u^\ve (x+\ve\omega)-u^\ve_\omega (x)$, from (\ref{eq3.1}) and (\ref{eq3.18}) by subtracting, we obtain 
the equation
\beq\label{eq3.19}
-\Div \l(a(x,\e^{-1}x
+\omega)\nabla w^\ve_\omega(x)\r)+w^\ve_\omega(x)=F^\e_{0,\omega}(x)+\Div F^\e_\omega(x),
\eeq
where
\beq\label{eq3.24}
\ds{
F^\e_{0,\omega}(x)=f(x+\ve\omega)-f(x),
}\atop\ds{
F^\e_\omega(x)=\l[a\l(x+\ve\omega,\e^{-1}x
+\omega\r)-a\l(x,\e^{-1}x
+\omega\r)\r]\nabla u^\ve_\omega (x).}
\eeq
We can write the energy inequality for  (\ref{eq3.19}) 
$$
\|w^\ve_\omega\|_{H^1(\RR^d)}\le c_0 (\|F^\e_{0,\omega}\|_{H^{-1}(\RR^d)}+\|F^\e_\omega\|_{L^2(\RR^d)}),\quad c_o=conct(\lambda).
$$
Hence, 
 using the property (\ref{sh2}) of the shift, the Lipschitz continuity of $a(x,y)$
and the energy estimate of the type (\ref{ener}) for the solution $u^\ve_\omega (x)$, we derive  firstly the inequality in $\rd$
\[
\|w^\ve_\omega\|^2_{H^1(\RR^d)}\le c \ve^2\|f\|^2_{L^2(\RR^d)} \quad \forall\omega\in Y,\quad c=const(d,\lambda,c_L).
\]
Integrating it over $\omega\in Y$ and recalling that 
$w^\ve_\omega (x)=u^\ve (x+\ve\omega)-u^\ve_\omega (x)$, we come to
\beq\label{eq3.20}
\int\limits_Y \|u^\ve (\cdot+\ve\omega)-u^\ve_\omega (\cdot)\|^2_{H^1(\RR^d)}
d\omega
\le  c \ve^2\|f\|^2_{L^2(\RR^d)}.
\eeq
 Comparing (\ref{eq3.20}) 
  with (\ref{eq3.5}) yields (\ref{eq3.5a}), by the triangle inequality.
  Lemma \ref{LemAv1} is proved.
  
  \bigskip
4$^\circ$ To prove Lemma \ref{LemAv2}, we apply elementary transformations in the left-hand side of (\ref{eq3.5a}).   We first change the  variable of integration in the integral over $\rd$, next change the order of integration, and finally use the convexity argument. Namely, 
$$
\int\limits_{Y}\int\limits_{\RR^d}(|u^\ve (x+\ve\omega)-u(x)-\e U(x,\f{x}{\ve}+\omega)|^2+| \ldots  |^2)
dx\,d\omega
$$
$$
=\int\limits_{Y}\int\limits_{\RR^d}(|u^\ve (x)-u(x-\ve\omega)-\e U(x-\ve\omega,\f{x}{\ve})|^2+| \ldots  |^2)
dx\,d\omega
$$
$$
=\int\limits_{\RR^d}\int\limits_{Y}|u^\ve (x)-u(x-\ve\omega)-\e U(x-\ve\omega,\f{x}{\ve})|^2d\omega\,dx +\int\limits_{\RR^d}\int\limits_{Y}| \ldots  |^2
d\omega\,dx
$$
$$
\ge\int\limits_{\RR^d}|\int\limits_{Y}\l(u^\ve (x)-u(x-\ve\omega)-\e U(x-\ve\omega,\f{x}{\ve})\r)d\omega|^2\,dx +\int\limits_{\RR^d}\int\limits_{Y}| \ldots  |^2
d\omega\,dx
$$
$$
=\int\limits_{\RR^d}|u^\ve (x)-\int\limits_{Y}u(x-\ve\omega)d\omega-\e \int\limits_{Y}U(x-\ve\omega,\f{x}{\ve})d\omega|^2\,dx +\int\limits_{\RR^d}|\int\limits_{Y} \ldots d\omega |^2
\,dx
$$
(for brevity,  we do not show  explicitly  transformations   in the second term with the gradient in (\ref{eq3.5a}), for they  are
quite clear and repeat those that are shown in the first term).
We see above the smoothed corrector  (\ref{eq3.5c}), that is
\[
K_\e(x)=\int\limits_{Y}U(x-\ve\omega,\f{x}{\ve})d\omega,
\] 
 and  Steklov's smoothing of the solution to the homogenized equation
\[
(S^\e u)(x)=\int\limits_{Y}u(x-\ve\omega)d\omega
\]
which can be replaced with the solution $u(x)$ itself, by the property 
(\ref{st1}) of 
the Steklov smoothing operator and the elliptic estimate (\ref{ells}).

In summary, the estimate (\ref{eq3.5b}) is verified.

\section{Proof of the main result}

 The main result is formulated below in
Theorem \ref{Th5.1}. We divide our proof of it into several steps.

1$^\circ$ We start with 
the notation that will simplify rather cumbersome formulas.

Let $N=(N^1,\ldots,N^d)$, where $N^j$ is the solution  to (\ref{eq2.1}).  
We  denote: 
\begin{equation}\label{3.1}
N_{\e,\o}(x):=N(x, \f{x}{\e}+\o),\quad
 U^\e_\o(x):=N_{\e,\o}(x)\cdot\nab u(x). 
 \end{equation} 
 Then the estimate (\ref{eq3.5}) takes the form
\beq\label{3.2} 
\int\limits_{Y} \|u^\ve_\omega (\cdot)-u(\cdot)-\ve  U^\e_\o(\cdot)\|^2_{H^1(\RR^d)}d\omega\le C\ve^2 \|f\|^2_{L^2(\RR^d)}.
\eeq
In particular,
\[
 \|u^\ve_\omega -u-\ve  U^\e_\o\|_{Y\times\RR^d}\le C\ve \|f\|,
\]
and our  immediate goal will be to investigate the $L^2$-form
\beq\label{3.3}
(u^\ve_\omega -u-\ve  U^\e_\o,h)_{Y\times\RR^d}, \quad h\in L^2(\rd).
\eeq
Here and in the sequel,  we use the simplified 
notation for the inner product and the norm  in the spaces $\ld$ and
$L^2(Y\times\RR^d)$
\beq\label{3.4}\ds{
\|\cdot\|=\|\cdot \|_{\ld},\quad
(\cdot \,,\cdot\,)=(\cdot \,,\cdot\,)_{\ld},}
\atop\ds{
\|\cdot\|_{Y\times\RR^d}=\|\cdot \|_{L^2(Y\times\RR^d)},\quad
(\cdot \,,\cdot\,)_{Y\times\RR^d}=(\cdot \,,\cdot\,)_{L^2(Y\times\RR^d)}
.}
\eeq

We recall some facts about  homogenization of the  equation adjoint to (\ref{eq3.1}). That is
 \begin{equation}\label{3.5}
 v^\e_\o \in H^1 (\R^d), \quad 
(A^\o_{\ve})^*v_\o^\ve+v^\e_\o=h,\quad h{\in} L^2(\R^d),\quad \o\in Y,
\end{equation}
(this problem 
for $\o=0$ appeared earlier as the problem (\ref{1s})).
It is associated with the homogenized problem (\ref{homs});
the corresponding first approximation is of the form
\begin{equation}\label{3.6}
v(x)+\varepsilon V^\e_\o(x),\quad \mbox{ where }\quad V^\e_\o(x)=
\tilde{N}_{\e,\o}(x)\cdot\nab v(x),\quad \tilde{N}_{\e,\o}(x)=
\tilde{N}(x,x/\e+\o),
\end{equation}
and the vector $\tilde{N}$ is  composed of the solutions to the  cell problem 
(\ref{cps}).
What is more, the following estimate (that is a counterpart of (\ref{eq3.5}) or (\ref{3.2}))  holds
\beq\label{3.7} 
\int\limits_{Y} \|v^\ve_\omega (\cdot)-v(\cdot)-\ve  V^\e_\o(\cdot)\|^2_{H^1(\RR^d)}d\omega\le c\ve^2 \|h\|^2_{L^2(\RR^d)},\quad c=const(d,\lambda, c_L),
\eeq
with its simple corollary
\beq\label{3.8} 
\int\limits_{Y} \|v^\ve_\omega (\cdot)-v(\cdot)\|^2_{L^2(\RR^d)}d\omega\le c\ve^2 \|h\|^2_{L^2(\RR^d)}, \quad c=const(d,\lambda, c_L).
\eeq
To derive (\ref{3.8}) from (\ref{3.7}) it suffices  to use the estimate of the type (\ref{estcor}) for the corrector $V^\e_\o$.

In the sequel, we will refer to the energy and elliptic estimates relating to (\ref{3.5}) and (\ref{homs})
 respectively, those are
\beq\label{3.9}
\|v^\ve_\o \|_{ H^1(\R^d)}\le c \|f\|\quad \forall\o\in Y, \quad 
c=cost(\lambda), 
\eeq
\beq\label{3.10}
\|v \|_{ H^2(\R^d)}\le c \|f\|, \quad 
c=cost(\lambda). 
\eeq

2$^\circ$ To investigate the $L^2$-form
(\ref{3.3}) we  insert $ u_\o^\e-u-\varepsilon U_\o^\e$ as a test function
into the integral identity
for the solution of the adjoint equation 
(\ref{3.5}), integrate it over $\o\in Y$  and make some transformations:
\[
(u^\ve_\omega -u-\ve  U^\e_\o,h)_{Y\times\RR^d}=
(u^\ve_\omega -u-\ve  U^\e_\o,((A^\o_{\ve})^*+1)v^\e_\o)_{Y\times\RR^d}
\]
\[
=((A^\o_\e+1)u^\e_\omega -(A^\o_\e+1)(u+\e  U^\e_\o),v^\e_\o)_{Y\times\RR^d}=
\l((A_0+1)u -(A^\o_\e+1)(u+\e  U^\e_\o),v^\e_\o\r)_{Y\times\RR^d}
\]
\beq\label{3.11}
=(A_0u-A^\o_\e(u+\e  U^\e_\o),v^\e_\o)_{Y\times\RR^d}-
\e( U^\e_\o,v^\e_\o)_{Y\times\RR^d}=:T_1-T_2.
\eeq

We study first the term $T_2$ in (\ref{3.11}). Since
\[
T_2:=\e( U^\e_\o,v^\e_\o)_{Y\times\RR^d}=\e( U^\e_\o,v^\e_\o-v)_{Y\times\RR^d}+\e( U^\e_\o,v)_{Y\times\RR^d},
\]
and
\[
( U^\e_\o,v)_{Y\times\RR^d}\stackrel{(\ref{3.1})}=(N_{\e,\o}\cdot\nab u,v)_{Y\times\RR^d}=(\int_Y N(x,x/\e+\o)d\,\o\cdot\nab u,v)=0
\]
because $\langle N(x,\cdot)\rangle=0$, we deduce, by the H\"{o}lder inequality, that
\[
|T_2|\le \e\|U^\e_\o\|_{Y\times\RR^d}\|v^\e_\o-v\|_{Y\times\RR^d},
\]
where
\[
\|U^\e_\o\|_{Y\times\RR^d}\stackrel{(\ref{estcor})}\le c \|f\|,\quad 
\|v^\e_\o-v\|_{Y\times\RR^d}\stackrel{(\ref{3.8})}\le c \e \|h\|.
\]
Thus, we conclude that
\beq\label{3.12}
T_2\cong 0. 
\eeq
Here and in the sequel, we  use the sign $\cong$ to denote any equality modulo terms $T$ having the following estimate
\[
|T|\le c\e^2\|f\|\,\|h\|,\quad c=const(d,\lambda,c_L);
\]
and such terms $T$ will be called  inessential. 

We proceed now to the more difficult term  $T_1$ in (\ref{3.11}).
We need the relations similar to (\ref{eq3.6})--(\ref{solG}) where the shifted functions like
\[
N_{\e,\o}=N(x,\f{x}{\e}+\o),\,g^j_{\e,\o}=g^j(x,\f{x}{\e}+\o),\,
G^j_{\e,\o}=G^j(x,\f{x}{\e}+\o),\,a^\e_\o=a(x,\f{x}{\e}+\o)
\]
 are involved. We don`t  formulate  here these  "shifted"\ 
relations, but refer to them by numbers (corresponding to their counterparts with $\o=0$) endowed with the index $\o$. For example,
there holds the representation
\beq\label{3.13}
\ds{
a^\e_\o\nab (u+\e U^\e_\o) - a^0\nab u\stackrel{(\ref{eq3.6})_\o}=
g^j_{\e,\o}\f{\pa u}{\pa x_j}+\e a^\e_\o\nab_x(N_{\e,\o}\cdot \nab u))
}\atop\ds{
\stackrel{(\ref{eq3.10})_\o}=\e\Div\l(G^j_{\e,\o}\f{\pa u}{\pa x_j}\r)-
\e\Div_x\l(G^j_{\e,\o}\f{\pa u}{\pa x_j}\r)+\e a^\e_\o\nab_x(N_{\e,\o}\cdot \nab u)
,} 
\eeq
where, for brevity, we denote 
\beq\label{3.14}
\ds{
\Div_x\l(G^j_{\e,\o}\f{\pa u}{\pa x_j}\r)=\f{\pa u(x)}{\pa x_j}
\Div_x\,G^j(x,y+\o)+G^j(x,y+\o)\nab \f{\pa u(x)}{\pa x_j},\quad
y=\f{x}{\e},}\atop\ds{
\nab_x(N_{\e,\o}\cdot \nab u)=\f{\pa u(x)}{\pa x_j}\nab_x N^j(x,y+\o)+
 N^j(x,y+\o)\nab \f{\pa u(x)}{\pa x_j},\quad y=\f{x}{\e}
.} 
\eeq
Therefore,
\[
T_1:=(A_0u-A^\o_\e(u+\e  U^\e_\o),v^\e_\o)_{Y\times\RR^d}=-
(a^\e_\o\nab (u+\e U^\e_\o) - a^0\nab u,\nab v^\e_\o)_{Y\times\RR^d}
\]
\beq\label{3.15}
\stackrel{(\ref{3.13})+(\ref{solG})_\o}=
\l(\e\Div_x\l(G^j_{\e,\o}\f{\pa u}{\pa x_j}\r)-\e a^\e_\o\nab_x(N_{\e,\o}\cdot \nab u),\nab v^\e_\o\r)_{Y\times\RR^d} 
\eeq
\[
\cong\l(\e\Div_x\l(G^j_{\e,\o}\f{\pa u}{\pa x_j}\r)-\e a^\e_\o\nab_x(N_{\e,\o}\cdot \nab u),\nab( v+\e V^\e_\o)\r)_{Y\times\RR^d}.
\]
We have just deleted the term
\[
T:=\l(\e\Div_x\l(G^j_{\e,\o}\f{\pa u}{\pa x_j}\r)-\e a^\e_\o\nab_x(N_{\e,\o}\cdot \nab u),\nab (v^\e_\o -v-\e V^\e_\o)\r)_{Y\times\RR^d}
\]
which is inessential. We show this,  using the H\"{o}lder inequality:
\[
|T|\le\e\,\|\Div_x\l(G^j_{\e,\o}\f{\pa u}{\pa x_j}\r)-\e a^\e_\o\nab_x(N_{\e,\o}\cdot \nab u)\|_{Y\times\RR^d}
\|\nab (v^\e_\o -v-\e V^\e_\o)\|_{Y\times\RR^d}\cong 0,
\]
because
\beq\label{3.16}
\|\nab (v^\e_\o -v-\e V^\e_\o)\|_{Y\times\RR^d}\stackrel{(\ref{3.7})}
\le c\e\|h\|,\quad c=const(d,\lambda,c_L),
\eeq
and
\beq\label{3.17}
\|\Div_x\l(G^j_{\e,\o}\f{\pa u}{\pa x_j}\r)-\e a^\e_\o\nab_x(N_{\e,\o}\cdot \nab u)\|_{Y\times\RR^d}
\le c\|f\|,\quad c=const(d,\lambda,c_L).
\eeq
 The latter can be derived 
 by  arguments used in the proof of (\ref{estcor}), if we take into account the structure of the functions (\ref{3.14}) involved in (\ref{3.17}), the elliptic estimate (\ref{ells}),
  and the properties of the oscillating factors in (\ref{3.17}) listed in Lemma \ref{locper} and Lemma \ref{locper1}.
  
  Returning to (\ref{3.15}), we continue to study the term $T_1$. First of all, we restore the vector $g^j_{\e,\o}$ in it, by using (\ref{eq3.10})$_\o$ and (\ref{solG})$_\o$:
\[
\l(\e\Div_x\l(G^j_{\e,\o}\f{\pa u}{\pa x_j}\r),\nab (v+\e V^\e_\o)\r)_{Y\times\RR^d}
\] 
\[
\stackrel{(\ref{eq3.10})_\o}=
\l(\e\Div\l(G^j_{\e,\o}\f{\pa u}{\pa x_j}\r),\nab (v+\e V^\e_\o)\r)_{Y\times\RR^d}
-\l(g^j_{\e,\o}\f{\pa u}{\pa x_j},\nab (v+\e V^\e_\o)\r)_{Y\times\RR^d}
\]
\[
\stackrel{(\ref{solG})_\o}=-\l(g^j_{\e,\o}\f{\pa u}{\pa x_j},\nab (v+\e V^\e_\o)\r)_{Y\times\RR^d}.
\]
Thus,
\beq\label{3.18}
T_1\stackrel{(\ref{3.15})}\cong
-(g^j_{\e,\o}\f{\pa u}{\pa x_j},\nab (v+\e V^\e_\o))_{Y\times\RR^d}
-\e
\l( a^\e_\o\nab_x(N_{\e,\o}\cdot \nab u),\nab (v+\e V^\e_\o)\r)_{Y\times\RR^d} =:I+II.
\eeq
Engaging the equality
\beq\label{3.19}
\nab (v+\e V^\e_\o)
\stackrel{(\ref{3.6})}=(\nab_y\tilde{N}^k_{\e,\o}+e^k)\f{\pa v}{\pa x_k}+\e\nab_x(\tilde{N}_{\e,\o}\cdot \nab v)
\eeq
with
\[
\nab_y\tilde{N}^k_{\e,\o}:=\nab_y\tilde{N}^k(x,y+\o),\quad y=\f{x}{\e},
\]
and $\nab_x(\tilde{N}_{\e,\o}\cdot \nab v)$ quite similar to  (\ref{3.14})$_2$,
we have the representation
\beq\label{3.20}
I=-\l(g^j_{\e,\o}\cdot (\nab_y\tilde{N}^k_{\e,\o}+e^k) \f{\pa u}{\pa x_j},\f{\pa v}{\pa x_k}\r)_{Y\times\RR^d}-
\e\l(g^j_{\e,\o} \f{\pa u}{\pa x_j},\nab_x(\tilde{N}_{\e,\o}\cdot 
\nab v)\r)_{Y\times\RR^d}.
\eeq
The oscillating vector in the first summand of 
(\ref{3.20}) has zero mean value with respect to  $\o$. In fact,
\[
\langle g^j_{\e,\o}\cdot (\nab_y\tilde{N}^k_{\e,\o}+e^k)\rangle_\o=
\langle g^j_{\e,\o}\cdot \nab_y\tilde{N}^k_{\e,\o}\rangle_\o+
\langle g^j_{\e,\o}\rangle_\o\cdot e^k
\]
\[
=\langle g^j(x,y+\o)\cdot \nab_y\tilde{N}^k(x,y+\o)\rangle_\o+
\langle g^j(x,y+\o)\rangle_\o\cdot e^k, \quad y=x/\e,
\]
if we explicate the notation introduced for brevity.
Since $\nab_y\tilde{N}^k(x,y+\o)=\nab_\o\tilde{N}^k(x,y+\o)$,  the both mean values in the last sum are equal to zero, in view of (\ref{solg})$_\o$. Consequently,
\[
\l(g^j_{\e,\o}\cdot (\nab_y\tilde{N}^k_{\e,\o}+e^k) \f{\pa u}{\pa x_j},\f{\pa v}{\pa x_k}\r)_{Y\times\RR^d}=
\l(\langle g^j_{\e,\o}\cdot (\nab_y\tilde{N}^k_{\e,\o}+e^k)\rangle_\o \f{\pa u}{\pa x_j},\f{\pa v}{\pa x_k}\r)=0,
\]
and (\ref{3.20}) yields
\beq\label{3.21}
I\cong-
\e\l(g^j_{\e,\o} \f{\pa u}{\pa x_j},\nab_x(\tilde{N}_{\e,\o}\cdot 
\nab v)\r)_{Y\times\RR^d}.
\eeq

To study the term $II$ in (\ref{3.18}), we insert the representation
(\ref{3.19}) in it. Then
\beq\label{3.22}\ds{
II=-
\e\l(\nab_x(N_{\e,\o}\cdot \nab u),(a^\e_\o)^*\nab (v+\e V^\e_\o)\r)_{Y\times\RR^d}
}\atop\ds{
\cong
-
\e\l(\nab_x(N_{\e,\o}\cdot \nab u),(a^\e_\o)^*(\nab_y\tilde{N}^k_{\e,\o}+e^k)\f{\pa v}{\pa x_k}\r)_{Y\times\RR^d},
}\eeq
where the inessential term is deleted, that is
\[
T:=\e^2\l(\nab_x(N_{\e,\o}\cdot \nab u),(a^\e_\o)^*\nab_x(\tilde{N}_{\e,\o}\cdot \nab v)\r)_{Y\times\RR^d}\cong
0.\]
To prove the last "approximate" equality, we use the H\"{o}lder inequality
and
 similar arguments, as in the proof of (\ref{estcor}), if we take into account the structure of the functions $\nab_x(N_{\e,\o}\cdot \nab u)$ and $(a^\e_\o)^*\nab_x(\tilde{N}_{\e,\o}\cdot \nab v)$ involved in $T$ (see, e.g., (\ref{3.14})$_2$).
 
We introduce the counterpart of the vector $g^j$ (see (\ref{eq3.80})) for the adjoint equation, that is
\[
\tilde{g}^j(x,y)= a^*(x,y)(e^j+\nabla_y \tilde{N}^j(x,y))-\langle a^*(x,\cdot)(e^j+\nabla_y \tilde{N}^j(x,\cdot))\rangle_Y,\]
or, rewriten in view of   (\ref{hom}), 
\[
\tilde{g}^j(x,y)= a^*(x,y)(e^j+\nabla_y \tilde{N}^j(x,y))- (a^0(x))^*e^j,
\]
whence
\beq\label{3.23}
(a^\e_\o)^*(\nab_y\tilde{N}^k_{\e,\o}+e^k)=\tilde{g}^k_{\e,\o}+(a^0)^*e^k.
\eeq
Inserting (\ref{3.23}) in (\ref{3.22}) yields
\beq\label{3.24}
II\cong-
\e\l(\nab_x(N_{\e,\o}\cdot \nab u),\tilde{g}^k_{\e,\o}\f{\pa v}{\pa x_k}\r)_{Y\times\RR^d},
\eeq
where we have deleted the term
\[
\e\l(\nab_x(N_{\e,\o}\cdot \nab u),(a^0)^*e^k\f{\pa v}{\pa x_k}\r)_{Y\times\RR^d}=\e\l(\nab_x(N_{\e,\o}\cdot \nab u),(a^0)^*\nab v\r)_{Y\times\RR^d}
\]
\[
=-\e\l(N_{\e,\o}\cdot \nab u,\Div_x((a^0)^*\nab v\r)_{Y\times\RR^d}
\stackrel{(\ref{3.1})}=
-\e\l(\int_Y N(x,x/\e+\o)\,d\,\o\cdot \nab u,\Div_x((a^0)^*\nab v\r)=0,
\]
 because $\langle  N(x,\cdot)\rangle=0 $.
 
 Collecting (\ref{3.18}), (\ref{3.21}) and (\ref{3.24}) together, we have
\beq\label{3.25}
 T_1\cong
-\e\l(g^j_{\e,\o} \f{\pa u}{\pa x_j},\nab_x(\tilde{N}_{\e,\o}\cdot 
\nab v)\r)_{Y\times\RR^d}-
\e\l(\nab_x(N_{\e,\o}\cdot \nab u),\tilde{g}^k_{\e,\o}\f{\pa v}{\pa x_k}\r)_{Y\times\RR^d},
\eeq
 where, according to the notation (\ref{3.14})$_2$,
 \[
 \nab_x(N_{\e,\o}\cdot \nab u)=\f{\pa u(x)}{\pa x_j}\nab_x N^j(x,y+\o)+
 N^j(x,y+\o)\nab \f{\pa u(x)}{\pa x_j},\quad y=\f{x}{\e},
 \]
 and similarly
  \[
 \nab_x(\tilde{N}_{\e,\o}\cdot \nab v)=\f{\pa v(x)}{\pa x_k}\nab_x \tilde{N}^k(x,y+\o)+
 \tilde{N}^k(x,y+\o)\nab \f{\pa v(x)}{\pa x_k},\quad y=\f{x}{\e}.
 \]
 Therefore, putting the oscillating factors close to each other and integrating by parts in (\ref{3.25}), we obtain
 \[
 T_1\cong
-\e\l(\f{\pa }{\pa x_k}\Div\langle\tilde{N}^k g^j\rangle \f{\pa u}{\pa x_j},v)\r)+
\e\l(\f{\pa }{\pa x_k}\langle\nab_x \tilde{N}^k\cdot g^j\rangle \f{\pa u}{\pa x_j}, v\r)
 \]
 \[
-\e\l(u,\f{\pa }{\pa x_j}\Div\langle\tilde{g}^k N^j\rangle \f{\pa v}{\pa x_k})\r)+
\e\l(u,\f{\pa }{\pa x_j}\langle\nab_x N^j\cdot  \tilde{g}^k\rangle \f{\pa v}{\pa x_k}\r),
 \]
 or shortly
 \beq\label{3.26}
 T_1\cong
-\e\l(L_3u-L_2u,v)\r)-
\e\l( u,\tilde{L}_3v-\tilde{L}_2v\r),
\eeq
where $L_3$, $L_2$, $\tilde{L}_3$, $\tilde{L}_2$ are differential operators of order three or two (which is indicated in index) with coefficients depending only on the \ "slow"\ variable $x$. Namely,
 \beq\label{3.27}
 \ds{L_3=D_k D_m c^{jk}_m(x)D_j,\quad \tilde{L}_3=
 D_j D_m  \tilde{c}^{kj}_m(x)D_k,
 }\atop\ds{
c^{jk}_m(x) =\langle g^j_m(x,\cdot) \tilde{N}^k (x,\cdot) \rangle,
\quad\tilde{c}^{kj}_m(x)=\langle \tilde{g}^k_m(x,\cdot) N^j(x,\cdot)\rangle,
}
 \eeq
  \beq\label{3.28}
 \ds{L_2=D_k c^{jk}(x)D_j,\quad \tilde{L}_2=
 D_j  \tilde{c}^{kj}(x)D_k,
 }\atop\ds{
c^{jk}(x) =\langle g^j(x,\cdot) \cdot \nab_x\tilde{N}^k (x,\cdot) \rangle,
\quad\tilde{c}^{kj}(x)=\langle \tilde{g}^k(x,\cdot)\cdot \nab_x N^j(x,\cdot)\rangle.
}
 \eeq
 Here we use the notation $D_j=\f{\pa}{\pa x_j}$, $j=1,\ldots,d$.
 
From  (\ref{3.11}), (\ref{3.12}), (\ref{3.26}), we get
  \beq\label{3.29}
 (u^\ve_\omega -u-\ve  U^\e_\o,h)_{Y\times\RR^d}\cong
 -\e((L_3-L_2) u,v)-\e(u,(\tilde{L}_3-\tilde{L}_2)v),
 \eeq
 where the operators $L_3$, $L_2$, $\tilde{L}_3$, $\tilde{L}_2$ are defined in (\ref{3.27}), (\ref{3.28}).
 
 3$^\circ$ Now we slightly change the form (\ref{3.3}), replacing in it the function $u^\e_\o$ with the shifted solution to the original problem
 that is
 \[
 u^\e(x+\e\o)=:S^\e_\o u^\e,\quad \o\in Y.
 \]
 We come to the form
 \[
(S^\e_\o u^\ve -u-\ve  U^\e_\o,h)_{Y\times\RR^d}, \quad h\in L^2(\rd),
\]
which we write as a sum
\beq\label{3.30}
(S^\e_\o u^\ve -u-\ve  U^\e_\o,h)_{Y\times\RR^d}=
 (S^\e_\o u^\ve -u^\e_\o,h)_{Y\times\RR^d}+
 (u^\e_\o -u-\ve  U^\e_\o,h)_{Y\times\RR^d}.
\eeq
We begin to estimate the first summand,
 addressing the solution of the adjoint equation (\ref{3.5}), as in the chain of equalities (\ref{3.11}):
\[
(S^\e_\o u^\ve -u^\e_\o,h)_{Y\times\RR^d}=(S^\e_\o u^\ve -u^\e_\o,((A_\e^\o)^*+1)v^\e_\o)_{Y\times\RR^d}
\] 
\[
=((A_\e^\o+1)(S^\e_\o u^\ve -u^\e_\o),v^\e_\o)_{Y\times\RR^d}.
\]
Recalling that the function 
$w^\ve_\omega (x)=u^\ve (x+\ve\omega)-u^\ve_\omega (x)=S^\e_\o u^\ve (x) -u^\e_\o (x)$ satisfies the equation (\ref{eq3.19}), which we write here in the form
\[
(A_\e^\o+1)(S^\e_\o u^\e -u^\e_\o)=S^\e_\o f -f+\Div\, F^\e_\o,
\]
where
\beq\label{3.31}
\ds{
F^\e_\omega(x)
=\l[a\l(x+\ve\omega,\e^{-1}x
+\omega\r)-a\l(x,\e^{-1}x
+\omega\r)\r]\nabla u^\ve_\omega (x)=S^\e_\o\l(a^\e(x)-\tilde{a}^\e(x)\r)\nabla u^\ve_\omega (x),
}\atop\ds{
a^\e(x)=a(x,\e^{-1}x),\quad \tilde{a}^\e(x)=a(x-\e\o,\e^{-1}x).
}
\eeq
Hence,
\beq\label{3.32}
\ds{
(S^\e_\o u^\ve -u^\e_\o,h)_{Y\times\RR^d}=
(S^\e_\o f -f,v^\e_\o)_{Y\times\RR^d}-(S^\e_\o\l(a^\e-\tilde{a}^\e\r)\nabla u^\e_\omega,\nab v^\e_\o)_{Y\times\RR^d}
}\atop\ds{
=:T^{(1)}-T^{(2)}.}
\eeq
Engaging the approximation (\ref{3.6}), we write the sum
\[
T^{(1)}=(S^\e_\o f -f,v^\e_\o-v-\e V^\e_\o)_{Y\times\RR^d}+
(S^\e_\o f -f,v+\e V^\e_\o)_{Y\times\RR^d}
\]
\[
\cong (S^\e_\o f -f,v+\e V^\e_\o)_{Y\times\RR^d},
\]
where we have deleted the  inessential term, due to the estimates (\ref{3.7}) and (\ref{sh2}).
Further transformations give
\beq\label{3.33}
\ds{
T^{(1)}\cong
(S^\e_\o f -f,v)_{Y\times\RR^d}+\e(S^\e_\o f -f,V^\e_\o)_{Y\times\RR^d}
}\atop\ds{
=(f,S^\e v -v)-\e(f,V^\e_\o)_{Y\times\RR^d}
+\e(S^\e_\o f ,V^\e_\o)_{Y\times\RR^d},}
\eeq
where there  emerges the Steklov smoothing operator
\beq\label{3.34}
S^\e v(x)=\int_Y  v(x-\e\o)\,d\omega.
\eeq
Note that
\[
\|S^\e v -v\|\stackrel{(\ref{st2})}\le c\e^2\|\nab^2v\|\stackrel{(\ref{3.10})}\le C\e^2
\|h\|
\]
and
\[
(f,V^\e_\o)_{Y\times\RR^d}\stackrel{(\ref{3.6})}=(f,\int_Y N(x,\f{x}{\e}+\o)\,d\omega\cdot \nab v)=0
\] because  $\langle N(x,\cdot)\rangle=0$.
Therefore, (\ref{3.33}) yields
\beq\label{3.35}
T^{(1)}\cong \e(S^\e_\o f ,V^\e_\o)_{Y\times\RR^d}.
\eeq

Let us proceed to the term $T^{(2)}$ in (\ref{3.32}). 
Since
\beq\label{3.36}\ds{
S^\e_\o\l(a^\e-\tilde{a}^\e\r)\stackrel{(\ref{3.31})}=
a(x+\e\omega,y
+\omega)-a(x,y 
+\omega)
}\atop\ds{
=\e\int_0^1
\nab_x a(x+t\e\o,y+\o)\cdot\o\,d\,t,\quad y=\frac x\ve,
}
\eeq
then
\[
T^{(2)}:=(S^\e_\o\l(a^\e-\tilde{a}^\e\r)\nabla u^\e_\omega,\nab v^\e_\o)_{Y\times\RR^d}
\]
\[
\cong (S^\e_\o\l(a^\e-\tilde{a}^\e\r)\nabla (u+\e U^\e_\omega),\nab (v+
\e V^\e_\o))_{Y\times\RR^d}\]
\[
\cong
\l(S^\e_\o\l(a^\e-\tilde{a}^\e\r)(\nab_y{N}^j_{\e,\o}+e^j)\f{\pa u}{\pa x_j}
,(\nab_y\tilde{N}^k_{\e,\o}+e^k)\f{\pa v}{\pa x_k}\r)_{Y\times\RR^d},
\]
where at each step some inessential terms are dropped away.
Among the deleted terms, there are
\[
(S^\e_\o\l(a^\e-\tilde{a}^\e\r)\nabla (u^\e_\omega-u-\e U^\e_\omega),\nab v^\e_\o)_{Y\times\RR^d},\quad
(S^\e_\o\l(a^\e-\tilde{a}^\e\r)\nabla (u+\e U^\e_\omega),\nab( v^\e_\o-v-
\e V^\e_\o))_{Y\times\RR^d},
\]
\[
\e(S^\e_\o\l(a^\e-\tilde{a}^\e\r)\nab_x(N_{\e,\o}\cdot \nab u),\nab( v^\e_\o-v-
\e V^\e_\o))_{Y\times\RR^d}
\]
and others.
 To show that these terms are inessential, we use quite standard for this paper arguments repeated not once. 
Except for (\ref{3.36}), we refer here to the 
estimates (\ref{3.2}), (\ref{3.7}), (\ref{ells}), (\ref{3.9}), (\ref{3.10}), the properties of the cell problems solutions, the H\"{o}lder inequality and, certainly, Lemma \ref{LemShift}.

Taking into account (\ref{3.36}), we come to the representation
\beq\label{3.37}\ds{
T^{(2)}
\cong
\e\l(\int_0^1
\nab_x a(x+t\e\o,y+\o)|_{y=x/\e}\cdot\o\,d\,t
(\nab_y{N}^j_{\e,\o}+e^j)\f{\pa u}{\pa x_j}
,(\nab_y\tilde{N}^k_{\e,\o}+e^k)\f{\pa v}{\pa x_k}\r)_{Y\times\RR^d}
}\atop\ds{
=\e\l(
\hat{c}^{jk}_\e\f{\pa u}{\pa x_j}
,\f{\pa v}{\pa x_k}\r)=:(
M_\e u,v),\quad 
M_\e=-\Div\, \hat{c}^{jk}_\e\nab,
}
\eeq
where the 
matrix $\{\hat{c}^{jk}_\e(x)\}_{jk}$ is obtained from $\nab_x a(x,y)$ through the procedure of double averaging     with respect to the both variables $x$ and $y$. Namely,
\beq\label{3.38}
\hat{c}^{jk}_\e(x)=\int_Y (\nab_y\tilde{N}^k_{\e,\o}+e^k)\cdot\l(\int_0^1
\nab_x a(x+t\e\o,y+\o)|_{y=x/\e}\cdot\o\,d\,t\r)
(\nab_y{N}^j_{\e,\o}+e^j)\,d\,\o
\eeq
with
\[
\nab_y{N}^j_{\e,\o}=\nab_y{N}^j(x,y+\o),\quad
\nab_y\tilde{N}^k_{\e,\o}=\nab_y\tilde{N}^k(x,y+\o),\quad y=\f{x}{\e}.
\]
We see that the  integration in (\ref{3.38}) does not consume fully  the parameter $\e$.

4$^\circ$ 
Gathering the estimates (\ref{3.30}), (\ref{3.29}), (\ref{3.32}), (\ref{3.35}), (\ref{3.37}), we arrive at
\beq\label{3.39}
(S^\e_\o u^\ve -u-\ve  U^\e_\o,h)_{Y\times\RR^d}\cong
 -\e((L_3-L_2) u,v)-\e(u,(\tilde{L}_3-\tilde{L}_2)v)-\e(M_\e u,v)+ \e(S^\e_\o f ,V^\e_\o)_{Y\times\RR^d}.
\eeq

We can replace (modulo inessential terms) the function $h(x)$  with  its shifting $S^\e_\o h(x)=
 h(x+\e\o)$  in the left-hand side form in (\ref{3.39}), by the property
 (\ref{sh2}) of the shift operator  and Lemma \ref{LemAv1}. Then
 \[
(S^\e_\o u^\ve -u-\ve  U^\e_\o,h)_{Y\times\RR^d}\cong
(S^\e_\o u^\ve -u-\ve  U^\e_\o,S^\e_\o h)_{Y\times\RR^d}
 \]
 \[
 =\int\limits_{Y\times \RR^d} (S^\e_\o u^\ve -u-\ve  U^\e_\o)S^\e_\o h\,dx\,d\o \stackrel{(\ref{3.1})}=
 \int\limits_{Y}\l(\int\limits_{\RR^d} ( u^\ve(x) -u(x-\e\o)-\ve  U(x-\e\o,\f{x}{\e}) h(x)
 \,dx\r)d\o
 \]
  \[
 =
 \int\limits_{\RR^d} \l(\int\limits_{Y}( u^\ve(x) -u(x-\e\o)-\ve  U(x-\e\o,\f{x}{\e}) h(x)
 \,d\o\r)\,dx
 \]
  \[
 =
 \int\limits_{\RR^d} \l(( u^\ve(x) -\int\limits_{Y}u(x-\e\o)\,d\o-\ve  \int\limits_{Y}U(x-\e\o,\f{x}{\e})\,d\o 
 \r)h(x)\,dx
 \]
 \[
\stackrel{(\ref{st})+(\ref{eq3.5c})+(\ref{3.1})} =(u^\ve-S^\e u-\e K_\e,h),
 \]
 where we have done the transformations inside the integral form similar to those used for the proof of Lemma \ref{LemAv2} (see the very end of the \S4). In view of (\ref{st2}), $S^\e u$ can be replaced with $u$ in the last form, and so we obtain
\beq\label{3.40}
(S^\e_\o u^\ve -u-\ve  U^\e_\o,h)_{Y\times\RR^d}\cong
(u^\ve- u-\e K_\e,h)=(u^\ve- u-\e \K_\e f,h)
\eeq
with $K_\e=\K_\e f$ defined in (\ref{eq3.5c}) and (\ref{eq3.5o}).

Similarly,
\beq\label{3.41}\ds{
(S^\e_\o f ,V^\e_\o)_{Y\times\RR^d}=
\int\limits_{Y\times \RR^d}  f(x+\e\o) V(x,\f{x}{\e}+\o)\,dx\,d\o
=\int\limits_{Y\times \RR^d}  f(x) V(x-\e\o,\f{x}{\e})\,dx\,d\o
}
\atop\ds{
=
\int\limits_{ \RR^d} f(x)
\l(\int\limits_{Y}V(x-\e\o,\f{x}{\e})\,d\o\r)\,dx=(f,\tilde{K}_\e)=
(f,\tilde{\K}_\e h)
}
\eeq
with the smoothed corrector 
\beq\label{3.42}
\tilde{K}_\e(x)=\int\limits_{Y}V(x-\e\o,\f{x}{\e})\,d\o=
\int\limits_{Y}\tilde{N}(x-\e\o,\f{x}{\e})\cdot \nab v(x-\e\o)\,d\o=:
\tilde{\K}_\e h
\eeq
for the adjoint problem (\ref{1s}).

From (\ref{3.39})--(\ref{3.41}), it follows that
\beq\label{3.43}
(u^\ve- u-\e \K_\e f,h)\cong
-\e((L_3-L_2) u,v)-\e(u,(\tilde{L}_3-\tilde{L}_2)v)-\e(M_\e u,v)+ \e(f,\tilde{\K}_\e h).
\eeq
Introducing
\beq\label{3.44}\ds{
\L:=(A_0+1)^{-1}\l(L_3-L_2+(\tilde{L}_3-\tilde{L}_2)^*\r)(A_0+1)^{-1},
}\atop\ds{
\mathcal{M}_\e:=(A_0+1)^{-1}M_\e(A_0+1)^{-1},
}
\eeq
we rewrite (\ref{3.43}) as follows
\[
((A_\e+1)^{-1}f- (A_0+1)^{-1}f-\e (\K_\e +\tilde{\K}_\e^*-\mathcal{L}-\mathcal{M}_\e)f,h)\cong
0,
\]
which means exactly the estimate sought
\beq\label{3.45}
\|(A_\e+1)^{-1}f- (A_0+1)^{-1}f-\e (\K_\e +\tilde{\K}_\e^*-\mathcal{L}-\mathcal{M}_\e)f\|\le c\e^2\|f\|,\quad c=const(d,\lambda,c_L),
\eeq
if the meaning of the symbol $\cong$ is taken into account. Obviously,
the estimate (\ref{3.45})
is equivalent to (\ref{eq1.9}) with the correcting term (\ref{eq1.10}).

We have proved
\begin{teor}\label{Th5.1} 
Under  assumptions  (\ref{eq1.2})--(\ref{eq1.4}), 
the resolvent $(A_\e+1)^{-1}$ is approximated with the sum  
$(A_0+1)^{-1}+\e (\K_\e +\tilde{\K}_\e^*-\mathcal{L}-\mathcal{M}_\e)$
so that the estimate (\ref{3.45}) holds true. The terms $\K_\e$, $\tilde{\K}_\e$, $\mathcal{L}$, $\mathcal{M}_\e$
 of the corrector are defined in (\ref{eq3.5c}) and (\ref{eq3.5o}), in
 (\ref{3.42}), in (\ref{3.27}), (\ref{3.28})  and (\ref{3.44})$_1$, in
 (\ref{3.37}), (\ref{3.38})  and (\ref{3.44})$_2$, respectively.
 \end{teor}
 
\noindent \textbf{Remark 6.2}. 
Coefficients of the operators defined in (\ref{3.27}) and (\ref{3.28}) are actually  calculated in terms of only the solutions ${N}^j$, $\tilde{N}^k$ to the cell problems
(\ref{eq2.1}), (\ref{cps}), respectively, and their gradients either.
No additional cell problems are needed. (Note that this is also valid for other components of the correcting operator in (\ref{3.45}) which is seen directly from their definitions.)
In fact, 
since $g^j_m=g^j\cdot e^m$, we have
\[
c^{jk}_m(x)\stackrel{(\ref{3.27})} =\langle g^j_m(x,\cdot) \tilde{N}^k (x,\cdot) \rangle
\stackrel{(\ref{eq3.80})}=\langle \tilde{N}^k(x,\cdot) ( a(x,\cdot) (\nab{N}^j(x,\cdot) +e^j)-a^0(x)e^j)\cdot e^m \rangle\]
\[=\langle \tilde{N}^k(x,\cdot) a(x,\cdot)(\nab{N}^j(x,\cdot)+e^j)\cdot e^m\rangle- \langle \tilde{N}^k(x,\cdot)\rangle a^0(x)e^je^m\]
\[\stackrel{(\ref{cps})}
=\langle \tilde{N}^k(x,\cdot) a(x,\cdot)(\nab{N}^j(x,\cdot)+e^j)\rangle
\cdot e^m.
\]
Similar expressions can be found for the coefficients of  the other operators from (\ref{3.27}) and (\ref{3.28}).

\bigskip
\noindent \textbf{Remark 6.3}. 
Let coefficients of the operator $A_\e$ in (\ref{0.1}) oscillate over two different groups of variables with different small periods $\e$ and $\delta=\delta(\e)$. We assume that $\delta/\e$ tends to zero as $\e$ tends to zero. It is known that the limit problem is obtained through reiterated homogenization procedure and corresponds to an elliptic equation with constant coefficients. The difference for resolvents of the original and the limit operators is estimated in operator $\ld$-norm; this estimate is of order $\max\{\e,\delta/\e\}$ (see, e.g., \cite{PT},  \cite{PaT15} or \cite{UMN}). The resolvent approximations of higher order can be found  by method we demonstrate here. This may be a subject for a subsequent paper.

\end{document}